\documentclass[a4paper,11pt]{amsart}

\usepackage{graphicx}
\usepackage{url}
\usepackage{amssymb}
\usepackage{amsthm}
\usepackage{amsmath}
\usepackage{stmaryrd}
\usepackage{paralist}
\usepackage{hyperref}
\usepackage{anysize}

\newcommand\fl{\mathcal{L}}
\renewcommand\mod{\;\mathsf{mod}\;}
\renewcommand\L{\mathsf{L}}
\renewcommand\dim{\mathsf{dim}\,}
\newcommand\R{\mathbb{R}}
\newcommand\K{\mathsf{K}}
\newcommand\Z{{\mathbb{Z}}}
\newcommand\Sind{\mathsf{ind_{SK}}\,}
\newcommand\KG{\mathsf{KG}}
\newcommand\nf{\mathcal{F}}
\newcommand\conv{\mathsf{conv}\,}
\newcommand\A{\mathcal{A}}
\renewcommand\l{\lambda}
\newcommand\tl{\ensuremath{\text{\boldmath{$\l$}}}}
\renewcommand\P{\mathcal{P}}
\renewcommand\H{\mathcal{H}}
\newcommand\Lskel{\Sigma_k}
\newcommand\fType[1]{\Lambda_k(#1)}
\newcommand\Skel[2]{\Sigma_{#1}({#2})}
\newcommand\Edim[1]{\mathsf{e\text{-}dim}({#1})}
\newcommand\simplex[1]{\Delta_{#1}}             
\newcommand\Psimplex{\simplex{n-1}^r}           

\newcommand\wedgeproduct{\leftslice}
\renewcommand\wp[2]{\mathcal{W}_{{#1},{#2}}}
\newcommand\WPsurf[2]{\mathcal{S}_{{#1},{#2}}}
\renewcommand\min{\mathsf{min}}
\renewcommand\max{\mathsf{max}}
\newcommand\eqset{I}

\newtheorem*{thm*}{Theorem}
\newtheorem{thm}{Theorem}[section]

\newtheorem{cor}[thm]{Corollary}
\newtheorem*{cor*}{Corollary}
\newtheorem{lem}[thm]{Lemma}
\newtheorem{prop}[thm]{Proposition}

\theoremstyle{definition}
\newtheorem{dfn}[thm]{Definition}
\newtheorem{example}[thm]{Example}

\newtheorem*{obs*}{Observation}

\theoremstyle{remark}
\newtheorem{rem}[thm]{Remark}

\title{Non-projectability of polytope skeleta}

\author{Thilo R\"{o}rig}
\address{Institut f\"{u}r Mathematik, MA 8-3, TU Berlin, 10623 Berlin, 
Germany}
\email{roerig@math.tu-berlin.de}

\author{Raman Sanyal}
\address{Department of Mathematics, UC Berkeley, 970 Evans Hall,
Berkeley, CA 94720}
\email{sanyal@math.berkeley.edu}

\thanks{\noindent Thilo R\"{o}rig was supported by DFG Research Group Polyhedral
Surfaces.  Raman Sanyal was supported by the Deutsche Forschungsgemeinschaft
within the research training group `Methods for Discrete Structures' (GRK
1408) and a Miller Research Fellowship.}

\keywords{dimensional ambiguity, projection of polytopes, topological
obstructions, products of polygons, products of simplices, wedge products, polyhedral surfaces}

\date{\today}

\begin{document}

\begin{abstract}\noindent 
    We investigate necessary conditions for the existence of projections of
    polytopes that preserve full $k$-skeleta.  More precisely, given the
    combinatorics of a polytope and the dimension $e$ of the target space, what
    are obstructions to the existence of a geometric realization of a polytope
    with the given combinatorial type such that a linear projection to
    $e$-space strictly preserves the $k$-skeleton.  Building on the work of
    Sanyal (2009), we develop a general framework to calculate obstructions to
    the existence of such realizations using topological combinatorics. Our
    obstructions take the form of graph colorings and linear integer programs.
    We focus on polytopes of product type and calculate the obstructions for
    products of polygons, products of simplices, and wedge products of
    polytopes. Our results show the limitations of constructions for the
    deformed products of polygons of Sanyal \& Ziegler (2009) and the wedge
    product surfaces of R\"{o}rig \& Ziegler (2009) and complement their results.
\end{abstract}

\maketitle

\section{Introduction}

According to Gr\"unbaum \cite[Ch.~12]{grue03}, a polytope $P$ is
\emph{dimensionally $k$-ambiguous} if the $k$-skeleton of $P$ is isomorphic to
that of a polytope $Q$ and $\dim Q \not= \dim P$.  So, not only is the
$k$-skeleton of such a polytope \emph{not} characteristic but, even worse, it
does not even give away the dimension in which to look for it. Unfortunately,
there is no effective way to decide when a polytope is dimensionally ambiguous
and even the list of known instances of such polytopes is rather short.  The
prime example of a dimensionally $\lfloor\frac{d-3}{2}\rfloor$-ambiguous
polytope is the $d$-simplex as is certified by the existence of
\emph{neighborly simplicial polytopes} such as the cyclic polytopes (cf.\
\cite{zie95}). However, in recent years two more families of polytopes joined
the list: the family of cubes via the existence of \emph{neighborly cubical
polytopes} \cite{jos00} and the family of products of even polygons in guise
of \emph{projected deformed products of polygons} \cite{zie04, sz07}. In both
cases, the construction principle (unified in \cite{sz07}) is to give a
special realization of the combinatorial type and to verify that a projection
to lower dimensions strictly preserves the skeleton in question.

The main motivation for this paper was to investigate the limitations of this
approach. To be more precise: What are necessary conditions for the existence
of a polytope $P \subset \R^d$ of a fixed combinatorial type and a projection
$\pi: \R^d \rightarrow \R^e$ such that $P$ and $\pi(P)$ have isomorphic
$k$-skeleta.

Building on technology developed in \cite{san07}, we devise tools that give
fairly good necessary conditions on the existence of such pairs $(P,\pi)$ in
terms of topological combinatorics. The main observation is that if $\pi: P
\rightarrow \pi(P)$ retains the $k$-skeleton for $k \ge 0$ then there is an
associated pair of spaces $(\partial\A,\|\Sigma_k\|)$ with $\|\Sigma_k\|
\hookrightarrow \partial\A$ where $\Sigma_k$ is a simplicial complex and
$\partial\A$ is a (polyhedral) sphere. The simplicial complex $\Sigma_k$ is
defined in terms of the combinatorics of $P$ whereas the dimension of the
sphere $\partial\A$ depends on $e$.  Thus, the existence of $(P,\pi)$ implies
that $\Sigma_k$ is embeddable into a sphere of a specific dimension.
Obstructing the embeddability of $\Sigma_k$ into a sphere of this 
dimension then impedes the existence of $(P,\pi)$. Drawing from methods of
topological combinatorics \cite{MatousekBZ:BU}, our obstructions take the form
of graph coloring problems and integer linear programs.

We focus on polytopes of \emph{product type} for which the factorization of
the skeleta allows us to replace the simplicial complex $\Sigma_k$ by somewhat
simpler subcomplexes. We apply the tools to the following three classes of
polytopes: 

\begin{asparaenum}

    \item[\it Products of polygons.] One curiosity left in connection with the
        deformed products of polygons of \cite{sz07} is that the general
        construction scheme fails for \emph{odd} polygons,
        i.e.\ polygons with an odd number of vertices.  With respect to the
        number of \emph{even and odd} polygons we prove necessary conditions
        on products of polygons to be dimensionally ambiguous via projection.
        Along the way, we obtain interesting byproducts.  For example, it is
        known, though apparently nowhere written up properly, that there is no
        realization of a product of two odd polygons such that a projection to
        the plane retains all vertices. As a teaser, we generalize this result
        to

        \noindent
        \emph{There is no realization of a product of $r$ odd polygons such
        that a projection to $r$-space retains all vertices.}

    \item[\it Products of simplices.] Products of simplices are ubiquitous in
        geometric and topological combinatorics. Most notable are their
        appearances in work on tropicial geometry and subdivisions
        \cite{santos05}, game theory and polynomial equations \cite{sturm02},
        and as building blocks for prodsimplicial complexes such as
        Hom-complexes \cite{pfeifle07}. It is known to both discrete geometers
        and topologists that no $d$-polytope is dimensionally $k$-ambiguous
        for $k \ge \lfloor \frac{d}{2} \rfloor$ (cf.~Theorem
        \ref{thm:vanKampen}).  Essentially, the reason is that the statement
        is already false for the $d$-simplex.  In
        Section~\ref{sec:ProductsOfSimplices}, we investigate obstructions to
        the projectability of products of simplices -- calculating these
        obstructions is intricately related to the coloring of \emph{Kneser
        graphs}. We generalize a result in~\cite{san07} that products of $r
        \ge d$ simplices of dimension $d$ cannot retain all vertices under
        projection to lower dimensions.  

      \item[\it Wedge products.] The properties of (combinatorial) products that
        we exploit for the calculation of the obstructions hold for more general
        polytope constructions, most notably the wedge product. The wedge
        product, introduced in \cite[Ch.~4]{Z109} (see also \cite{RZ09}), is a
        degeneration of the product that may be described purely
        combinatorially.  The interest for this class stems from the original
        context in which wedge products were introduced: The (straight)
        realization of (equivelar) polyhedral surfaces. The equivelar surfaces
        of type $\{r,2n\}$ are topological surfaces \emph{glued} exclusively
        from $r$-gons, $2n$ of which meet at every vertex. The
        discrete-geometric realization question now is to find a geometric
        embedding in which all the polygons are convex and flat.  In~\cite{Z109}
        it is shown that certain equivelar families of type $\{r,2n\}$ are
        naturally embedded into the wedge products~$\wp{r}{n-1}$ of $r$-gons and
        $(n-1)$-simplices.  Techniques similar to the deformed products
        (cf.~\cite{sz07}) allow for the realization of the subfamily $\{r,4\}$
        in euclidean $3$-space.  In Section~\ref{sec:WedgeProducts} we show
        (cf.\ Theorem~\ref{thm:ObstructionWedgeProductSurface}) that this is
        probably the only family that embeds into $3$-space via projection:

        \noindent
    \emph{For $r \ge 4$ and $n \ge 3$ there is no realization of the wedge
    product $\wp{r}{n-1}$ such that a projection to $4$-space retains the
    surface $\WPsurf{r}{2n}$.}

    \noindent
    Our methods do not yield an obstruction for $r=3$ in which case the
    surface is triangulated and the wedge product of triangle and
    $(n-1)$-simplex is a simplex.
\end{asparaenum}

{\bf Acknowledgments.} We would like to thank G\"{u}nter Ziegler for
stimulating discussions and comments on an earlier draft of this paper.

\section{Combinatorial types, projections, and obstructions}
\label{sec:projections}

In this section we develop a general framework for investigating the
projectability of skeleta or more general subcomplexes of the boundary of a
polytope. We briefly recap the necessary polytopal background and then proceed
to reduce polytopes to their combinatorial structure -- their combinatorial
types.  The benefit will be apparent in our results which state conditions
under which there is no polytope with specific combinatorial and geometric
qualities. 

Throughout a (convex) \emph{polytope} $P \subset \R^d$ is the convex hull of
finitely points $P = \conv \{ v_1,\dots, v_n \}$ and, equivalently, the
bounded intersection of finitely many halfspaces $P = \{ x \in \R^d : a_i\cdot
x \le b_i \text{ for all } i = 1,\dots,m \}$. In both representations, we assume
that the collection of \emph{vertices} $v_1,\dots, v_n$ and of
\emph{facet-defining inequalities} $a_i \cdot x \le b_i$ is irredundant, that
is, no vertex or linear inequality can be omitted. It is customary to write
the system of linear inequalities succinctly as $Ax \le b$.  A hyperplane $H =
\{ x \in \R^d : c \cdot x = \delta \}$ is \emph{supporting} $P$ if $P
\subseteq H^- = \{ x \in \R^d: c \cdot x \le \delta \}$ and $F = P \cap H$ is
called a \emph{face} of $P$ -- the emptyset and $P$ are also faces of~$P$. In particular, it is clear that every $a_i \cdot
x \le b_i$ is a supporting hyperplane and the corresponding faces are called
\emph{facets}.  The dimension $\dim F$ of a face $F \subseteq P$ is the
dimension of its affine span. Vertices are faces of dimension $0$ and facets
are faces of dimension $\dim P - 1$.  We abbreviate the notions of
$k$-dimensional face and $d$-dimensional polytope with $k$-face and
$d$-polytope, respectively.

\begin{prop}[{\cite[Prop.~2.3]{zie95}}]\label{prop:face_reps}
    Let $P \subset \R^d$ be a polytope with vertex set $V \subset \R^d$ and
    facets $F_i$ defined by $a_i \cdot x \le b_i$ for $i = 1,\dots,m$. If $F \subseteq P$ is a
    face, then
    \begin{enumerate}
        \item $F = \conv(F \cap V)$ and
        \item $F = \{ x \in P : a_i \cdot x = b_i \text{ for all } i \in
            I_P(F) \}$ with
            $
                I_P(F) := \{ i \in [m] : F \subseteq F_i \}
            $
    \end{enumerate}
\end{prop}

The collection of faces $\fl(P)$ of a polytope $P$ ordered by inclusion is
called the \emph{face lattice} of $P$. The face lattice is a graded lattice
of rank $\dim P -1$ and it can be thought of as the combinatorial structure of
$P$. We call two polytopes \emph{combinatorially isomorphic} if $\fl(P) \cong
\fl(Q)$ as graded lattices. It follows from Proposition~\ref{prop:face_reps}
that the face lattice has two canonical representations.

\begin{cor} Let $P$ be a polytope with vertex set $V$ and facets indexed by
    $[m] = \{ 1, 2, \dots, m\}$. Then $\fl(P)$ is isomorphic to 
    \begin{enumerate}
        \item $\{ F \cap V : F \subseteq P  \text{ a face}
            \} \subseteq 2^V$ ordered by inclusion and \hfill {\rm(vertex description)}
        \item $\{ I_P(F) : F \subseteq P  \text{ a face }  \} \subseteq
            2^{[m]}$ ordered by reverse inclusion. \hfill {\rm(facet
            description)}
    \end{enumerate} 
\end{cor}

Our main results will be concerned with the non-existence of geometric
realizations of polytopes with given combinatorial features under projection.
In order to avoid cumbersome formulations, we wish to abstract from the
geometry of a polytope $P$. 

\begin{dfn}
    A graded lattice $\P$ is called a \emph{combinatorial type} of dimension
    $d$, or \emph{$d$-type} for short, if $\P \cong \fl(P)$ for some
    $d$-polytope $P$. 
\end{dfn}

We want to think about combinatorial types as polytopes stripped from their
geometric realization but we will nevertheless stick to our geometric
terminology and, for example, call $F \in \P$ a face of $\P$. Moreover, when
no confusion arises we use $P$ and $\P$ interchangeably. Identifying the
collection of facets of $\P$ with $F_1, \dots, F_m$, we write 
\[
    I_\P(F) = \{ i : F \subseteq F_i \} \subseteq [m]
\]
for the \emph{facet-incidences} of $\P$. The collection of all faces of $\P$
of dimension at most $k$ is the \emph{$k$-skeleton} of $\P$ and we call a $d$-type
$\P$ \emph{simple} if every $k$-face $F$ is contained in exactly $d-k$ facets.

\subsection{Geometry and topology of projections}
\label{sec:GenCase}

Let $P$ be a $d$-polytope and let $\pi: P \rightarrow \pi(P) \subseteq \R^e$ be
an affine projection. Throughout it is understood that $d \ge e$ and that
$\pi(P)$ is full-dimensional.
We want to find conditions under which $P$ and $\pi(P)$ have isomorphic
$k$-skeleta. The key concept for establishing such conditions is that
of faces strictly preserved under $\pi$.

\enlargethispage{3em}
\begin{dfn}[Preserved and strictly preserved faces; {\cite{san07,zie04}}]
    Let $P$ be a polytope, $F \subset P$ a proper face and $\pi: P \rightarrow
    \pi(P)$ a projection of polytopes. The face $F$ is \emph{preserved} under
    $\pi$ if 
    \begin{enumerate}
        \item[  i)] $G = \pi(F)$ is a proper face of $\pi(P)$ and
        \item[ ii)] $F$ and $G$ are combinatorially isomorphic.
    \end{enumerate}
    If, in addition,
    \begin{enumerate}
        \item[iii)] $\pi^{-1}(G)$ is equal to $F$
    \end{enumerate}
    then $F$ is \emph{strictly preserved}.
\end{dfn}

With the notion of strictly preserved faces at our disposal, the task of
deciding isomorphic $k$-skeleta of $P$ and $\pi(P)$ can be checked one face at
a time.

\begin{lem}\label{lem:iso_skel}
    Let $P$ be a polytope and let $\pi: P \rightarrow \pi(P)$ be a projection
    of polytopes. For \mbox{$0 \le k < \dim P$} the polytopes  $P$ and
    $\pi(P)$ have isomorphic $k$-skeleta if and only if every $k$-face of $P$
    is strictly preserved under~$\pi$.
\end{lem}
\begin{proof} 
    Assume that $P$ and $\pi(P)$ have isomorphic $k$-skeleta. We show by
    induction on $k$ that all preserved $(k+1)$-faces are then strictly
    preserved. 

    Since $f_l(P) = f_l(\pi(P))$ for $0 \le l \le k$, the $0$-skeleton is
    strictly preserved. If for $l \ge 1$ the $(l-1)$-skeleton is strictly
    preserved under projection, then the preimage of every $l$-face of $\pi(P)$
    is an $l$-face.  Indeed, let $\bar{F}$ be an $l$-face of~$\pi(P)$ and let $F
    = \pi^{-1}(\bar{F})$.  Then the map $\pi|_{F}: F \rightarrow \pi(F) =
    \bar{F}$ is a projection of polytopes that strictly preserves the
    $(l-1)$-skeleton of $F$. Thus the $(l-1)$-skeleton of $F$ is a subcomplex of
    an $(l-1)$-sphere. Hence $F$ is an $l$-face of $P$ and $\bar{F}$ is strictly
    preserved.

    Therefore all $k$-faces are strictly preserved since all $k$-faces are
    preserved and $P$ and $\pi(P)$ have isomorphic $(k-1)$-skeleta.
    
    Conversely, since every $i$-face for $i \le k$ is strictly preserved we have
    that the \mbox{$k$-skeleton} of $P$ is isomorphic to a subposet of the
    $k$-skeleton of $\pi(P)$. Assume that the inclusion is strict and let $H
    \subset P$ be a proper face of dimension greater than $k$ and $\pi(H)$ a
    $k$-face of $\pi(P)$. As a polytope, $H$ has a proper face $F$ of dimension
    $k$. But~$F$ is a $k$-face of~$P$ with $\pi(F) = \pi(H)$, since $\pi(F)
    \subseteq \pi(H)$ and $\dim \pi(F) = \dim \pi(H) = k$. Thus $F$ is not
    strictly preserved.
\end{proof}

In \cite{san07}, for every simple polytope $P$ a simplicial complex $\Sigma_0
= \Skel{0}{P}$ is defined in terms of the combinatorics of the vertices of
$P$. Furthermore, it is shown that if $\pi : P \rightarrow \pi(P)$ is a
projection strictly preserving the vertices, then $\Sigma_0$ is realized as a
subcomplex of a (simplicial) sphere whose dimension depends on $\dim \pi(P)$.
Theorem~\ref{thm:main} below is a generalization of this result for which we
separate the technical part in the following proposition.

\begin{prop} \label{prop:main}
    Let $P$ be a $d$-polytope on $m$ facets and let $\pi : P \rightarrow \pi(P)
    $ be a projection retaining all vertices of $P$. Then there is a polytope
    $\A = \A(P,\pi)$ of dimension $m-d-1+\dim\pi(P)$ with vertices $a_1,
    a_2,\dots, a_m$ such that the following holds: For every strictly preserved
    face $G \subset P$ the set
    \[
        \A_G := \conv \{ a_i : i \in [m] \setminus \eqset_P(G) \}
    \]
    is a simplex face of $\A$.
\end{prop}
\begin{proof}
    Let $e = \dim \pi(P)$. Fix a strictly preserved face $G$ and let $I =
    \eqset_P(G)$.
    Proposition 3.8 and
    Lemma 3.2 in \cite{san07} assert that there exists a polytopal Gale
    transform $\mathcal{G} = \{ g_1, g_2, \dots, g_m \} \subset \R^{d-e}$ with
    the property that $\mathcal{G}_I := \{ g_i : i \in I \}$ positively
    spans $\R^{d-e}$.
    Let $\A = \conv\{a_1,\dots,a_m\}$ be the $(m-d-1+e)$-dimensional polytope
    Gale-dual to $\mathcal{G}$. Gale duality implies that $\A_G = \conv \{ a_i :
    i \not\in I\}$ is a face of $\A$. Clearly, every set $\mathcal{G}_J$ with $I
    \subseteq J \subseteq [m]$ is positively spanning as well. So Gale duality
    implies that every subset of the vertices of $\A_G$ is also a face. Hence
    $\A_G$ is a simplex face of $\A$.
\end{proof}

We call the polytope $\A(P,\pi)$ the \emph{projection polytope}.  The
collection of strictly preserved faces induces the following simplicial
complex in the boundary of $\A$ that certifies the strict preservation.

\begin{dfn}
    Let $P$ be a polytope on $m$ facets and let $\pi: P \rightarrow \pi(P)$ be
    a projection of polytopes retaining all vertices.  We define  the
    \emph{strict projection complex} $\K(P,\pi) \subseteq
    2^{[m]}$ as the simplicial complex generated by the sets $\{ [m] \setminus
    I_P(G) : G \text{ strictly preserved under } \pi\}.$
\end{dfn}

We may now rephrase Proposition \ref{prop:main} as follows.

\begin{thm}\label{thm:main}
    Let $P$ be a $d$-polytope on $m$ facets and let $\pi : P \rightarrow
    \pi(P)$ be a projection strictly preserving all vertices. Then $\K(P,\pi)$
    is embedded in a (polytopal) sphere of dimension \mbox{$m-d-2 + \dim
    \pi(P)$}.  \qed 
\end{thm}

\begin{rem}
    The conditions of Theorem \ref{thm:main} can be weakened to the requirement
    that for each facet $F$ there is a strictly preserved vertex $v$ with $v
    \not\in F$. The proof relies on a slight variation of \cite[Proposition
    3.8]{san07} which verifies that the set $\mathcal{G}$ is indeed a polytopal
    Gale transform.
\end{rem}

As we are primarily interested in the preservation of full skeleta of a given
dimension we introduce the following complex of a combinatorial type.

\begin{dfn}
    Let $\P$ be a combinatorial type of dimension $d$ on $m$ facets.  For $-1
    \le k \le d$, the \emph{$k$-th coskeleton complex} is the simplicial
    complex 
    \[
        \Skel{k}{\P} = \{ \tau \subseteq [m] : \tau \cap I_\P(G) = \emptyset
        \text{ for some $k$-face } G \in \P \} \subseteq 2^{[m]}.
    \]
\end{dfn}

The maximal faces of $\Skel{k}{\P}$ are in bijection with the $k$-faces
of $\P$ under the correspondence $G \mapsto [m] \setminus I_\P(G)$. The
connection to $\K(P,\pi)$ is the following.

\begin{obs*}\label{obs:tower}
    If $\pi : P \rightarrow \pi(P)$ is a projection retaining the $k$-skeleton
    then
    \[
        \{\emptyset\} = \Skel{-1}{P} \ \subset\  \Skel{0}{P} \ \subset\
        \Skel{1}{P} \ \subset\  \cdots \ \subset\  \Skel{k}{P} \ \subset\
        \K(P,\pi)
    \]
    is an increasing sequence of subcomplexes.
\end{obs*}

As every $k$-face is contained in at least $d-k$ facets, the dimension of
$\Skel{k}{\P}$ is at most $m + k - d - 1$. If $\P$ is a simple $d$-type, then
$\Skel{k}{\P}$ is pure of this dimension.  In \cite{san07}, $\Skel{0}{\P}$ was
defined for simple $d$-types in terms of the \emph{complement complex} of the
boundary complex of the dual of~$\P$. Here, we abandon the restriction to
simple polytopes. 

Every simplicial complex can be embedded in a sphere of some dimension.
We will be interested in the smallest dimension of such a sphere.

\begin{dfn}[Embeddability dimension] 
    Let $\K \subseteq 2^{[m]}$ be a simplicial complex on $m$ vertices. 
    The \emph{embeddability dimension} $\Edim{\K}$ is the smallest integer $d$ such
    that $\|\K\|$ may be embedded into the $d$-sphere, i.e.\ $\|\K\|$ is
    homeomorphic to a closed subset of the $d$-sphere.
\end{dfn}

Theorem \ref{thm:main} can be read as an upper bound on the embeddability
dimension of the strict projection complex $\K(P,\pi)$. However, $\K(P,\pi)$
heavily depends on the geometry of the projection and hence a priori our
knowledge about $\K(P,\pi)$ is rather limited.  The virtue of the coskeleton
complex is that it is a subcomplex of $\K(P,\pi)$ defined entirely in terms of
the combinatorics of~$P$.

\begin{cor}\label{cor:proj_obstruct}
    Let $\P$ be a $d$-type on $m$ facets and, for $0 \le k < d$, let $\Lskel =
    \Lskel(\P)$ be the $k$-th coskeleton complex of $\P$. If
    \[
        e \ <\ \Edim{\Sigma_k} + d - m + 2  
    \]
    then there is no realization of $\P$ such that a projection to $\R^e$
    retains the $k$-skeleton.
\end{cor}
\begin{proof}
    By contradiction, assume that $P$ is a realization of $\P$ and $\pi: P
    \rightarrow \pi(P)$ is a projection retaining the $k$-skeleton with $\dim
    \pi(P)  = e < \Edim{\Sigma_k} + d - m + 2$. By Theorem \ref{thm:main}, the
    above observation, and the fact that the embeddability dimension is
    monotone along subcomplexes, the complex $\Sigma_k$ is realized in a
    sphere of dimension
    \[
        \Edim{\Sigma_k} \le m - d - 2 + e < \Edim{\Sigma_k}.
    \]
\end{proof}

The following, well-known fact bounds the embeddability dimension of a
simplicial complex in terms of its dimension.

\begin{prop}[{\cite[Thm.\ 11.1.8, Ex.\ 4.8.25]{grue03}}]
    \label{prop:edim_triv_bnds}
    Let $\K$ be a simplicial complex of dimension $\dim \K = \ell$. Then
    \[
        \ell \ \le\ \Edim{\K} \ \le\ 2\ell + 1.
    \]
\end{prop}

It is instructive to consider the statement of
Corollary~\ref{cor:proj_obstruct} in the extreme cases of
Proposition~\ref{prop:edim_triv_bnds}.  If the $(m+k-d-1)$-dimensional complex 
$\Edim{\Sigma_k}$ attains the lower bound then
Corollary~\ref{cor:proj_obstruct} implies that the dimension of the target
space has to be at least $e \ge k+1$. This is reassuring as the projection
embeds $\Skel{k}{\P}$ into a sphere of dimension $e-1$.  Now, suppose that
$\Edim{\Sigma_k}$ attains the upper bound and that $\P$ is a simple type. Then
$\dim \Skel{k}{\P} = m - (d-k) - 1$ and the $k$-skeleton is not projectable to
$e$-space if $e < m-d+2k+1$. This is the linear Van~Kampen--Flores result:

\begin{thm}[{\cite[Thm.~2]{grue65}}]\label{thm:vanKampen}
    Let $\P$ be a $d$-type and let $0 \le k \le \lfloor \frac{d-2}{2}
    \rfloor$. If 
    \[
        e \ \le \ 2k + 1
    \] 
    then there is no realization of $\P$ such that a projection to $e$-space
    retains the $k$-skeleton.
\end{thm}

\subsection{Cotype complexes of products}
\label{sec:ProductTypes}

For our purposes we need better bounds than provided by
Proposition~\ref{prop:edim_triv_bnds} and so we need more sophisticated
techniques to determine or at least bound the embeddability dimension
$\Edim{\Sigma_k}$. In this and the next section we introduce two notions that
approximate the coskeleton complex as well as the embeddability dimensions and
allow us to calculate bounds.

For the cases in which we want to apply Corollary~\ref{cor:proj_obstruct}, the
combinatorial types under consideration are products or, at least, closely
related (cf. Section~\ref{sec:WedgeProducts}). Let $P \subset
\R^d$ and $P^\prime \subset \R^{d^\prime}$ be two polytopes of combinatorial types $\P$ resp.\
$\P^\prime$. The product of $P$ and $P^\prime$ is the polytope $P \times P^\prime = \conv \{ (p,p^\prime) :
p \in P, p^\prime \in P^\prime\}$.  Combinatorially we define
\[
    \P \times \P^\prime
    := \fl(P \times P^\prime) 
    = \{ (F,F^\prime) :
    F \in \fl(P), F^\prime \in \fl(P^\prime) \text{ such that  } F = \emptyset
    \text{ iff } F^\prime = \emptyset \}.
\]
Note that this product of combinatorial types differs from the usual direct product
of lattices inasmuch as every face of the product $\P\times \P^\prime$ is a
product of non-empty faces of $\P$ and $\P^\prime$. In particular, we have $\dim
(F,F^\prime) = \dim F \times F^\prime = \dim F + \dim F^\prime$ and the
facet incidences of the product are given by
\[
    I_{P \times P^\prime}(F \times F^\prime) = I_P(F) \uplus
    I_{P^\prime}(F^\prime).
\]
 The following definition distinguishes the faces of the product by their
 ``type''.

\begin{dfn}\label{dfn:FaceComplex}
    Let $\P = \P_1 \times \P_2 \times \cdots \times \P_r$ with $\P_i$ a
    $d_i$-type on $m_i$ facets for $i = 1,\dots,r$. For a fixed $0 \le k < d =
    d_1 + \cdots + d_r$  we call a composition $\tl = (\l_1,\dots,\l_r) \in
    \Z^r$ with $0 \le \l_i \le d_i$ and $\l_1 + \cdots + \l_r = k$ a
    \emph{face type} of dimension $k$. We denote by $\fType{\P}$ the
    collection of $k$-dimensional face types for $\P$.  For $\tl \in
    \fType{\P}$ we define the \emph{cotype complex} of type $\tl$ as the join
    of coskeleton complexes
    \[
        \Sigma_{\tl}(\P) := \Sigma_{\l_1}(\P_1) \ * \
                            \Sigma_{\l_2}(\P_2) \ * \
                            \cdots \ * \
                            \Sigma_{\l_r}(\P_r). 
    \]
\end{dfn}

It is clear from the definition of the product that every face of $\P$ belongs
to some type and this yields a partition of the coskeleton complex.

\begin{prop} \label{prop:ProductFaceTypes}
    Let $\P = \P_1 \times \P_2 \times \cdots \times \P_r$ and $0 \le k < \dim
    \P$. Then
    \[
        \Skel{k}{\P} = \bigcup_{\tl \in \fType{\P}} \Skel{\tl}{\P}.
    \]
    \\[-5ex]\qed
\end{prop}

The monotonicity of the embeddability dimension along subcomplexes yields our
first bound for the projectibility of products.

\begin{cor} \label{cor:projectabilityBound}
    Let $\P$ be a product and $0 \le k < \dim \P$. If there is a face type
    $\tl \in \fType{\P}$ such that
    \[
        e < \Edim{\Sigma_\tl} + d - m + 2 
    \]
    then there is no realization of $\P$ such that a projection to $\R^e$
    retains the $k$-skeleton. \qed
\end{cor}

\begin{example} 
    To illustrate the usefulness of the cotype complex, consider the following
    question: Is there a realization of $\P = \Delta_1 \times \Delta_2$, a
    prism over a triangle,  such that a projection to the plane preserves the
    three \emph{vertical} edges (see Figure \ref{fig:Desargue}).  The ad-hoc
    negation of the question is that by \emph{Desargues' Theorem}
    (cf.~\cite[Sect.~14.3]{cox89}) the three vertical edges in the prism meet
    in a common point (at infinity) and a linear projection retains this
    property. Using the developed machinery, we see that the assumed
    projection satisfies the conditions of Proposition \ref{prop:main} and
    the vertical edges correspond to the face type $\tl = (1,0)$. The cotype
    complex is also shown in Figure \ref{fig:Desargue}: it consists of three
    triangles that share a common edge. Corollary
    \ref{cor:projectabilityBound} implies that such a projection does not
    exist as $\Skel{(1,0)}{\P}$ is not planar, i.e.\ $\Edim{\Skel{(1,0)}{\P}} =
    3$. But $d=3$ and $m=5$ and hence Corollary~\ref{cor:projectabilityBound}
    yields the non-projectability because $\Edim{\Skel{(1,0)}{\P}} + 3 - 5 + 2 = 3>2=e$.

\begin{figure}[ht]
    \centering
         \includegraphics[width=.18\linewidth]{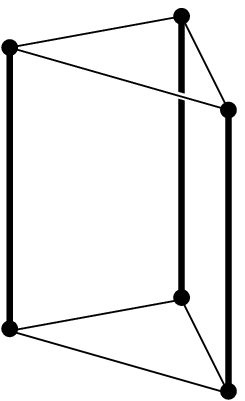}
         \hspace{.1\linewidth}
         \includegraphics[width=.18\linewidth]{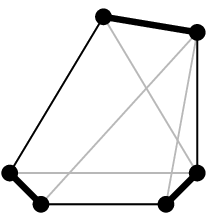}
         \hspace{.1\linewidth}
         \includegraphics[width=.18\linewidth]{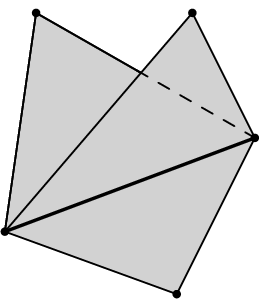}
         \caption{The triangular prism to the left with bold vertical edges.  An
           alleged projection in the middle with preserved vertical edges.  And
           the associated cotype complex to the right.}
        \label{fig:Desargue}
\end{figure}

\end{example}

\begin{rem}
    The definition of the cotype complex relies on properties of the product
    that are shared by other polytope constructions such as joins, direct sums,
    and wedge products (see Section \ref{sec:WedgeProducts}). The common
    generalization is that of a \emph{compound type} (cf.~\cite{raman-thesis})
    which is subject to further study.
\end{rem}

\subsection{Bounding the embeddability dimension}
\label{sec:SarkariaIndex}

In general it is hard to decide the embeddability of a complex $\K$ into some
$\R^e$. The following notions, taken and adapted from 
\cite{MatousekBZ:BU}, show that in fortunate cases bounds on $\Edim{\K}$
can be obtained combinatorially.

For a simplicial complex $\K \subseteq 2^{[m]}$ we denote by $\nf(\K)$ the set
of {\it minimal non-faces}, i.e.\ the inclusion-minimal sets in $2^{[m]}
\setminus \K$. The {\it Kneser graph} $\KG(\nf)$ on a set system $\nf
\subseteq 2^{[m]}$ has the elements of $\nf$ as vertices and $F,G \in \nf$
share an edge iff $F$ and $G$ are disjoint. Furthermore, for a graph $G$ we
denote by $\chi(G)$ the \emph{chromatic number} of $G$.

\begin{dfn}
    Let $\K$ be a simplicial complex on $m$ vertices and $\nf = \nf(\K)$ the
    collection of minimal non-faces. The {\it Sarkaria index} of $\K$ is 
    \[
        \Sind \K := m - \chi( \KG( \nf ) ) - 1.
    \]
\end{dfn}

\begin{thm}[Sarkaria's coloring/embedding theorem
    {\cite[Sect.~5.8]{MatousekBZ:BU}}] \label{thm:SarkariaColoringEmbedding}
    Let $\K$ be a simplicial complex. Then
    \[
        \Edim{\K} \ge \Sind \K.
    \]
\end{thm}

Every embedding of a simplicial complex $\K$ into a $d$-sphere gives rise to
a $\Z_2$-equivariant map of the \emph{deleted join} $\K^{*2}_\Delta$ to a
$d$-sphere. The \emph{$\Z_2$-index} of $\K^{*2}_\Delta$ is the smallest such
$d$ for which an equivariant map exists. In its original form
in~\cite{MatousekBZ:BU} the above theorem bounds from below the $\Z_2$-index
of $\K^{*2}_\Delta$ and thus also bounds from below the embeddability
dimension of $\K$.

The next observation reduces the calculation of the Sarkaria index of a
product to its factors.

\begin{prop}[{\cite[Prop.~3.10]{san07}}] \label{prop:IndexJoin}
    Let $\K$ and $\L$ be simplicial complexes. Then
    \[
        \Sind( \K * \L) = \Sind \K + \Sind \L + 1.
    \]
\end{prop}

Thus it follows directly from Definition~\ref{dfn:FaceComplex} that the Sarkaria
index of a cotype complex is determined by its factors.

\begin{cor}\label{cor:LinearIndex}
    Let $\P = \P_1 \times \P_2 \times \cdots \times \P_r$ and let $\tl \in
    \fType{\P}$. Then
    \[
        \Sind \Skel{\tl}{\P} = \sum_{i=1}^r \Sind \Skel{\l_i}{\P_i} + r - 1.
    \]
    \\[-2em]\qed
\end{cor}

We determine the exact embeddability dimensions and Sarkaria indices for two
coskeleton complexes of an arbitrary combinatorial type. The result depends only
on the number of facets.

\begin{prop} \label{prop:IndexSimplePolytopes}
    Let $\P$ be a $d$-type on $m$ facets. Then $\Sigma_d(\P) =
    \simplex{m-1}$ is homeomorphic to an $(m-1)$-ball and 
    \[
        m-1 = \Edim{\Skel{d}{\P}} = \Sind \Skel{d}{\P}.
    \]
    For the $(d-1)$-skeleton we have that
    $\Sigma_{d-1}(\P) = \partial \simplex{m-1} \cong S^{m-2}$ and 
    \[
        m-2 = \Edim{\Skel{d-1}{\P}} = \Sind \Skel{d-1}{\P}.
    \]
\end{prop}
\begin{proof}
    The first claim follows from the definition of the skeleton complex. Thus
    the embeddability dimensions are $m-1$ and $m-2$, respectively.  For the
    Sarkaria index we get in the former case that the Kneser graph of the
    minimal nonfaces of $\Sigma_d(\P)$ has no vertices, whereas in the latter
    case the graph has no edges. 
\end{proof}

In the special case that we have an $r$-fold product $\P^r = \P \times \P
\times \cdots \times \P$ of the same combinatorial
type $\P$, bounds on the embeddability dimension of $\Skel{k}{\P^r}$ can be
obtained by solving a \emph{knapsack-type} problem.

\begin{lem}\label{lem:knap}
    Let $\P$ be a $d$-type and let $r \ge 1$ and $0 \le k \le rd - 1$.
    For $i = 0, \dots, d$ set $s_i = \Sind \Skel{i}{\P}$ and let $s^*$ be the
    optimal value of the integer linear program
    \[
        \begin{array}{rr@{\ +\ }r@{\ +\ }r@{\ +\ }rl}
            \max & s_0\,\mu_0 & s_1\,\mu_1 & \cdots & s_d\,\mu_d \\
            \mathsf{s.t.}& 0\,\mu_0 & 1\,\mu_1 & \cdots & d\,\mu_d & =\ k \\
            & \mu_0 & \mu_1 & \cdots & \mu_d & =\ r \\
        \end{array}
    \]
    with $\mu_0,\dots,\mu_d \in \Z_{\ge 0}$. Then $\Edim{\Skel{k}{\P^r}} \ge
    s^*  + r - 1$.
\end{lem}
\begin{proof}
    To a face type $\tl \in \fType{\P^r}$ associate the non-negative numbers
    $(\mu_0, \mu_1,\dots,\mu_d)$ with 
    \[
        \mu_i \ =\ \#\{ j \in [r] : \lambda_j = i \}.
    \]
    They satisfy
    \[
        \begin{array}{r@{\ +\ }r@{\ +\ }r@{\ +\ }rl}
            0\,\mu_0 & 1\,\mu_1 & \cdots & d\,\mu_d & =\ k \text{ and} \\
                \mu_0 & \mu_1 & \cdots & \mu_d & =\ r \\
        \end{array}
    \]
    since $\tl$ is a partition of $k$ in $r$ parts and the Sarkaria index of
    $\Skel{\tl}{\P^r}$ is given by $\sum_i s_i\, \mu_i + r - 1$. Vice versa,
    every such non-negative collection of numbers $\mu_i$ that satisfies the
    conditions of the integer program gives rise to a valid face type.
\end{proof}

\section{Products of Polygons}
\label{sec:PolygonProducts}

Denote by $\P_m$ the combinatorial type of an $m$-gon, that is, a
$2$-dimensional combinatorial type on $m \ge 3$ facets labeled in cyclic
order. In this section we determine necessary conditions for the existence of
a realizations of $\P = \P_{m_1} \times \P_{m_2} \times \cdots \times
\P_{m_r}$, a product of polygons, that retain the $k$-skeleton under a
suitable projection. To that end, we need to determine (bounds on) the
embeddability dimension of $\Skel{k}{\P}$ for $0 \le k < 2r = \dim \P$. 

For a single polygon, Proposition~\ref{prop:IndexSimplePolytopes} leaves us to
determine the Sarkaria index for the $0$-th coskeleton complex of an $m$-gon.

\begin{lem} \label{lem:Sind_mgon} 
    Let $m \ge 3$ and $\P_m$ the combinatorial type of an $m$-gon. The
    Sarkaria index for the $0$-th coskeleton complex is given by
    \[
        \Sind \Sigma_0(\P_m) = \left\{
        \begin{array}{ll}
            m-3, & \text{ if $m$ is even, and } \\
            m-2, & \text{ if $m$ is odd. } \\
        \end{array}\right.
        \]
\end{lem}
\begin{proof}
    We show that the Kneser graph of minimal non-faces of $\Sigma_0(\P_m)$ has
    chromatic number $2$ and $1$, respectively.  For that let us
    determine the minimal non-faces of $\Skel{0}{\P_m}$: A subset $\sigma
    \subseteq [m]$ of the facets of $\P_m$ is a non-face of $\Sigma_0(\P_m)$ if
    and only if every vertex of $\P_m$ is incident to at least one facet $F_i$
    of $\P_m$ with $i \in \sigma$.  If a vertex of $\P_m$ is covered twice by
    $\sigma$ then every other minimal non-face intersects $\sigma$ and thus
    $\sigma$ is an isolated vertex in the Kneser graph. If $\sigma$ covers every
    vertex exactly once, then $[m] \setminus \sigma$ is again a minimal
    non-face.  It follows that for odd $m$ the Kneser graph consists of isolated
    vertices alone while for even $m$ there is exactly one edge.
\end{proof}

\begin{example}
    Let us consider $\Skel{0}{\P_5}$, the $0$-th coskeleton complex of the
    pentagon. The figure shows the triangles of the $0$-th coskeleton complex of
    the pentagon which form a M{\"o}biusstrip. Hence $\Skel{0}{\P_5}$ is not
    embeddable in the $2$-sphere.

    \begin{figure}[ht]
        \centering
        \includegraphics[width=.6\linewidth]{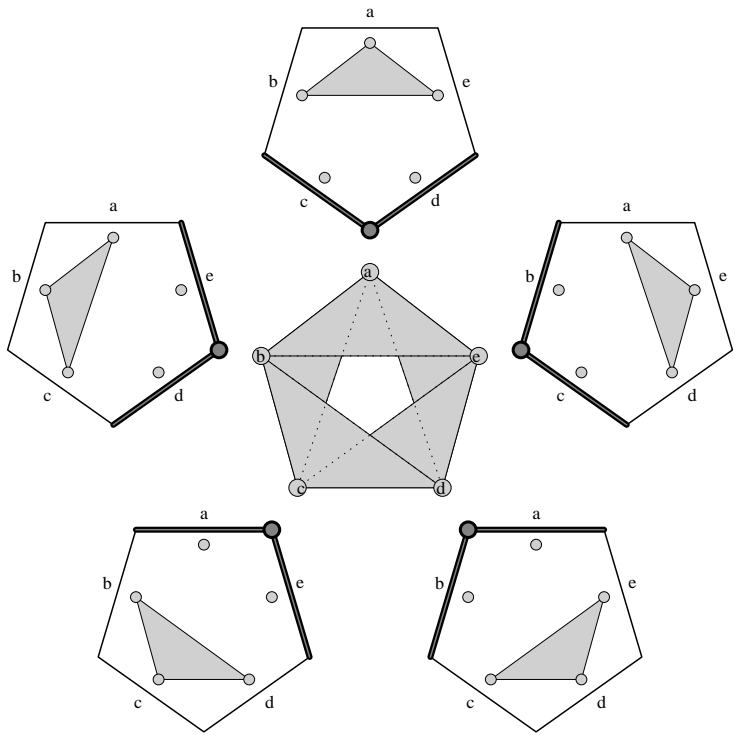}
        \vspace{.5cm}
        \caption{The five triangles of the $0$-th coskeleton complex of a
        pentagon fit together to form a M\"obiusstrip.}
        \label{fig:Moeb}
    \end{figure}

\end{example}

The example shows that the $0$-th coskeleton complex of an \emph{odd} polygon
has a certain twist to it that obstructs the embeddability into $m-3$
dimensional space.

Lemma~\ref{lem:Sind_mgon} implies that projectability bounds for products of
polygons arising from Corollary~\ref{cor:LinearIndex} will only depend on the
total number of facets and the number of odd and even polygons. Thus, it
suffices to consider the generic product
\[
    \P \ =\ \P_{\mathsf{even}}^{r_e} \ \times \ \P_{\mathsf{odd}}^{r_o},
\]
of $r_e$ \emph{even} and $r_o$ \emph{odd} polygons. We denote by $m$ the total
number of facets and by $r = r_e + r_o$ the number of factors. For a product
of polygons, we utilize the \emph{knapsack-type} integer program introduced in
Lemma~\ref{lem:knap}.

\begin{thm}\label{thm:EdimProductPolygons}
    Let $\P$ be a product of polygons with $r = r_o + r_e$ factors and $m$
    facets.  For $0 \le k \le 2r$ the embeddability dimension
    of the $k$-th coskeleton complex is bounded by
    \[
    \Edim{\Skel{k}{\P}} \ \ge \ m - 1 - r + \left\lfloor \frac{k}{2}
    \right\rfloor + \min\Big\{ 0, \left\lceil \frac{k}{2} \right\rceil - r_e
    \Big\}.
    \]
\end{thm}
\begin{proof}
    In the spirit of Lemma~\ref{lem:knap} consider the following
    integer linear program
    \[
    \newcommand\1{\phantom{1}}
    \newcommand\e{\mathsf{even}}
    \renewcommand\o{\mathsf{odd}}
    \begin{array}{llllll}
        \mathsf{min} & 2\,\mu_0^\e &+\ \1\,\mu_0^\o &+\
        \1\,\mu_1 \\
        \mathsf{s.t.} & & &\phantom{+}\
        \1\,\mu_1 &+\ 2\, \mu_2&=\ k\\
                      & \1\,\mu_0^\e &+\ \1\,\mu_0^\o &+\
        \1\,\mu_1 &+\ \1\, \mu_2&=\ r\\
                      & \1\,\mu_0^\e & & & &\le\ r_e\\
                      & &\phantom{+}\ \1\,\mu_0^\o & & &\le\ r_o\\
    \end{array}
    \]
    with $\mu_0^\mathsf{even}, \mu_0^\mathsf{odd},\mu_1,\mu_2 \in \Z_{\ge0}$.
    Every face type $\tl = (\lambda_1,\dots,\lambda_r)\in \Lambda_k(\P)$ gives
    rise to a feasible solution by the association
    \[
    \begin{array}{lll}
            \mu_2 &:=\ \# \{ i : \l_i = 2 \} & \text{ (polygons)} \\
            \mu_1 &:=\ \# \{ i : \l_i = 1 \} & \text{ (edges)} \\
            \mu_0^\mathsf{odd} &:=\ \# \{ i : \l_i = 0, 1 \le i \le r_o  \} & \text{ (odd
            vertices)} \\
            \mu_0^\mathsf{even} &:=\ \# \{ i : \l_i = 0, r_o < i \le r  \} & \text{ (even
            vertices)} \\
        \end{array}
    \]
    and, vice versa, every feasible solution yields a face type.  The integer
    program reduces to a problem in essentially two variables and the optimal
    solution is easily seen to be
    \[
        \mu^* = r - \left\lfloor \frac{k}{2} \right\rfloor + \max\Big\{ 0, r_e -
        \left\lceil\frac{k}{2}\right\rceil \Big\}.
        \]
    The result then follows from the fact that $\Edim{\Skel{k}{\P}} \ge
     m - 1 - \mu^*$.
\end{proof}

In order to put the above result in perspective, let us calculate upper bounds
on the embeddability dimension.

\begin{prop} \label{prop:prod_poly_upper}
    Let $\P$ be a product of polygons with $r = r_o + r_e$ factors and let $m$
    be the number of facets. For 
    $0 \le k < 2r$ the embeddability dimension is bounded by 
    \[
    \Edim{\Skel{k}{\P}} \le
    \left\{
    \begin{array}{ll}
        m - r - r_e - 1,& \text{ if } k = 0\\
        m - r - 1,& \text{ if } k = 1\\
        m - 2,& \text{ otherwise}.
    \end{array}
    \right.
    \]
\end{prop}
\begin{proof}
    \newcommand\bl{{\boldsymbol \ell}}
    \newcommand\Hskel{\hat{\Sigma}}
    Let $\P = \P_{m_e}^{r_e} \times \P_{m_o}^{r_o}$ and for
    $\ell = \min \{ k, 2 \}$ define
    \[
        \hat{\Sigma} \ =\ \Skel{\ell}{\P_{m_e}}^{*r_e} \ * \
        \Skel{\ell}{\P_{m_o}}^{*r_o}.
    \]
    We claim that $\hat{\Sigma}$ contains $\Sigma = \Skel{k}{\P}$ as a
    subcomplex. By construction, $\hat{\Sigma}$ and $\Sigma$ have identical
    vertex sets. For every admissible face type $\tl = (\l_1,\dots,\l_r) \in
    \fType{\P}$ we have $\l_i \le \ell$ for $i = 1,\dots,r$ and, by
    Observation \ref{obs:tower} and the relation of subcomplexes among joins,
    this shows $\Skel{\tl}{\P} \subseteq \Hskel$. Since $\Sigma$ is the union
    of all cotype complexes, this proves the claim.
    We will therefore bound the embeddability dimension of $\Hskel$ from above.

    For $\ell =  2$, we have by Proposition
    \ref{prop:IndexSimplePolytopes} that $\Skel{2}{\P_n} = \simplex{n-1}
    \hookrightarrow \partial\simplex{n}$ and thus $\Hskel$ embeds into the
    boundary of $\simplex{m_e}^{\oplus r_e} \oplus \simplex{m_o}^{\oplus r_o}$,
    a simplicial sphere of dimension $r_e m_e + r_o m_o - 1 = m - 1$.

    For $\ell = 1$, we again make use of  Proposition
    \ref{prop:IndexSimplePolytopes} to get $\Skel{1}{\P_n} =
    \partial\simplex{n-1}$ and therefore $\Hskel \hookrightarrow
    \partial(\simplex{m_e-1}^{\oplus r_e} \oplus \simplex{m_o-1}^{\oplus r_o})$,
    which is a simplicial sphere of dimension $r_e(m_e-1) + r_o(m_o-1) - 1 = m
    - r - 1$.

    For $\ell = 0$, the $0$-th coskeleton complex of $\P_n$ may be embedded into
    the boundary of an $(n-1)$-simplex. However, for even $n = 2t$ we can do
    better: Consider the $(n - 2)$-dimensional polytope $Q_t = \simplex{t-1}
    \oplus \simplex{t-1}$ and the mapping from the vertices of
    $\Skel{0}{\P_n}$
    that maps the $i$-th vertex to the $\lfloor\frac{i}{2}\rfloor$-th vertex of
    the first summand if $i$ is even and of the second otherwise. We claim
    that this gives an embedding. Every vertex $v$ of $\P_n$ is the
    intersection of an odd and an even edge. Thus the corresponding facet $[n]
    \setminus \eqset(v)$ is the disjoint union of $t-1$ odd and $t-1$ even
    vertices. These sets correspond to facets of $Q_t$. Thus $\Skel{0}{\P} =
    \Hskel$ embeds into the boundary of  $Q_t^{\oplus r_e} \oplus
    \simplex{m_o - 1}^{\oplus r_o}$ with $t = \frac{m_e}{2}$.
\end{proof}

Combining the bounds on the embeddability dimensions of the coskeleton
complexes of Theorem~\ref{thm:EdimProductPolygons} with
Corollary~\ref{cor:proj_obstruct} we obtain the following obstructions to
projectability of products of polygons.

\begin{thm}
    \label{thm:ProjectionProdPolygons}
    Let $r = r_o + r_e$ and  $0 \le k < 2r$. There is no realization
    of a product of $r_o$ odd and $r_e$ even polygons such that 
    a projection to $e$-dimensional space strictly preserves the $k$-skeleton
    if
    \[
    e < r + 1 + \left\lfloor \frac{k}{2} \right\rfloor +
    \min\Big\{0,\left\lceil\frac{k}{2}\right\rceil-r_e\Big\}.
    \]
    \\[-3em]\qed
\end{thm}

In~\cite{sz07}, $e$-dimensional polytopes with the
$\left\lfloor\tfrac{e-2}{2}\right\rfloor$-skeleton of the $r$-fold product of
\emph{even} polygons are projections of a suitable products of even polygons. The
following corollary shows that this construction technique does not generalize
to products of odd polygons.

\begin{cor}
    \label{cor:ProjectionOddPolygons}
    There is no realization of an $r_o$-fold product of \emph{odd} polygons
    such that the $k$-skeleton is strictly preserved under projection to
    $\R^e$ if
    \[
        e < r_o + 1 + \left\lfloor\frac{k}{2}\right\rfloor.
    \]
\end{cor}

In the special case of $r_o=2$ and $k=0$ the result reduces to the
well-known fact that a product of two odd $m$-gons does not project to an
$m^2$-gon.

Another case of interest is $k = \lfloor\tfrac{e}{2}\rfloor - 1$. 
In case such a realization and projection exists, the resulting polytope is
called neighborly, in analogy to the simplicial neighborly polytopes.

\begin{cor}
    \label{cor:ProjectionMoreEvenPolygons}
    Let $r = r_e + r_o$ and $e \ge 1$.  If 
    \[
    \left\{
    \begin{array}{rcll}
        \left\lceil\frac{3e-2}{4}\right\rceil &<& r 
        &\text{ for $r_e < \left\lfloor\frac{e}{4}\right\rfloor$,} \\[.5em]
        \left\lfloor\frac{e}{2}\right\rfloor &<& r_o 
        &\text{ for $r_e \ge \left\lfloor\frac{e}{4}\right\rfloor$}
    \end{array}\right.
    \]
    then there is no realization of a product of $r_e$ even and $r_o$ odd
    polytopes such that a projection to $e$-space is neighborly.
\end{cor}

Paraphrasing the situation for products of odd polygons, the result puts an
upper bound of $\lceil\tfrac{3e-2}{4}\rceil$ on the number of odd polygons for
a ``neighborly'' projection to $\R^e$.

\section{Products of simplices}
\label{sec:ProductsOfSimplices}

In this section we investigate obstructions to skeleta-preserving projections
of products of simplices.  Appealing to the results from Section
\ref{sec:projections}, we bound the embeddability dimension for the respective
coskeleton complexes. We denote by $\simplex{n-1} = 2^{[n]}$ the combinatorial
type of an $(n-1)$-simplex.  The key to determining the embeddability
dimension and the Sarkaria index of $\Skel{k}{\simplex{n-1}}$ is the following
observation.

\begin{obs*} 
    For $n \ge 1$ and $0 \le k \le n - 1$ the $k$-th coskeleton complex
    $\Skel{k}{\simplex{n-1}}$ of the $(n-1)$-simplex is isomorphic to the
    $k$-skeleton of $\simplex{n-1}$.
\end{obs*}

Thus $\Skel{k}{\simplex{n-1}}$ is a \emph{well known} complex and the
calculation of the Sarkaria index involves the \emph{classical} Kneser graphs
$\KG_{n,\ell} = \KG \tbinom{[n]}{\ell}$ for $0 \le \ell \le n$, that is, the
Kneser graphs on the collection of $\ell$-sets of an $n$-set. Their chromatic
numbers are a celebrated result in topological combinatorics.

\begin{thm}[Lov\'{a}sz~\cite{lov78}]\label{thm:lov_kneser}
    For $1 \le \ell \le n$ the chromatic number of $\KG_{n,\ell}$ is given by
    \[
        \chi(\KG_{n,\ell}) = 
        \begin{cases}
            n - 2\ell + 2& \text{if } \ell \le \tfrac{n+1}{2} \\
            1 & \text{otherwise.}
        \end{cases}
  \]
\end{thm}

This result immediately implies the Sarkaria index of the 
$k$-th skeleton complex $\Skel{k}{\simplex{n-1}}$. 

\begin{lem} \label{lem:IndexSimplex}
    For $n \ge 2$ and $0 \le k \le n-1$ the Sarkaria index of the $k$-th
    coskeleton complex $\Sigma_k = \Skel{k} {\simplex{n-1}}$ of the
    $(n-1)$-simplex is 
    \[
        \Sind \Sigma_k =
        \left\{
        \begin{array}{rl@{\ }cl}
            2k+1, & \text{ if } & 0              &\le\ k \ \le \ \tfrac{n-3}{2},  \\
            n-2,  & \text{ if } & \tfrac{n-3}{2} &< \ k \ \le \ n-2, \\
            n-1,  & \text{ if } & \multicolumn{2}{r}{ k \ = \ n-1.}
      \end{array}\right.
    \]
\end{lem}
\begin{proof}
    By the above observation, we have $\KG(\nf(\Sigma_k)) = \KG_{n,k+2}$.  The
    first two cases follow directly from Theorem~\ref{thm:lov_kneser}. The last
    case follows from Proposition~\ref{prop:IndexSimplePolytopes}.
\end{proof}

In combination with Proposition \ref{prop:edim_triv_bnds} we obtain the
following corollary.

\begin{cor}\label{lem:IndexSimplexUpper}
    Let $\Sigma_k = \Skel{k}{\simplex{n-1}}$ be its $k$-th skeleton complex of
    an $(n-1)$-simplex for $n \ge 2$.  Then the embeddability dimension
    satisfies
    \[
        \Edim{\Sigma_k} =
        \left\{
            \begin{array}{ll@{\ }cl}
            2 k+1,& \text{ if }& 0              &\le\  k \ \le\ \tfrac{n-3}{2}, \\
            n-2,  & \text{ if }& \tfrac{n-3}{2} &<\  k \ \le\ n-2,\\
            n-1,  & \multicolumn{3}{l}{\text{ otherwise.}}
        \end{array}
        \right.
    \]
    \qed
\end{cor}

In the following we denote by
\[
    \Psimplex = \underbrace{
    \simplex{n-1} \times
    \simplex{n-1} \times
    \cdots \times
    \simplex{n-1}}_{r}
\]
an $r$-fold product of $(n-1)$-simplices.

\begin{thm} \label{thm:IndexProdSimplices}
    Let $n \ge 2$, $r \ge 1$ and $0 \le k < r(n-1)$. The embeddability
    dimension of the $k$-th coskeleton complex $\Sigma_k = \Skel{k}{\Psimplex}$
    of the product of simplices $\Psimplex$ is bounded from below by
    \[
        \Edim{\Sigma_k} \ge
        \left\{
        \begin{array}{l@{\text{ if }}ccccc}
            2r + 2k - 1, & 
                0 & \le & k & \le & r \lfloor\tfrac{n-3}{2}\rfloor\\
                \tfrac{1}{2} rn + k - 1, &
                r\lfloor\tfrac{n-3}{2}\rfloor & < & k & \le &
                r\lfloor\tfrac{n-2}{2}\rfloor\\ 
            r(n-1)+\alpha-1, & 
                r\lfloor\tfrac{n-2}{2} \rfloor & < & k & < & r(n-1)
        \end{array}
        \right.
    \]
    and 
    \[
        \alpha = \left\lfloor \frac{ k - r \lfloor\frac{n-2}{2}\rfloor
        }{\lfloor\frac{n+1}{2}\rfloor}
        \right\rfloor.
    \]
\end{thm}

\begin{proof}
    We use the knowledge gained from Lemma \ref{lem:IndexSimplex} to set up
    the integer linear program as in Lemma~\ref{lem:knap}. Set $c =
    \lfloor\frac{n-3}{2}\rfloor$ and let  $0 \le k < r(n-1)$. The program is
    \[
        \begin{array}{rll}
            \max & \sum_{j=0}^{c} (2j+1)\, \mu_j + 
            (n-2) \sum_{j = c+1}^{n-2} \mu_j + (n-1) \mu_{n-1}\\
            \mathsf{s.t.} & \phantom{0\,}\mu_0 + \phantom{1\,}\mu_1 + \cdots
            +\phantom{(n-1)\,}\mu_{n-1}&= r\\
            & 0\,\mu_0 + 1\,\mu_1 + \cdots + (n-1)\mu_{n-1} &= k\\
        \end{array}
    \]
    and subject to the condition that the $\mu_i$ are non-negative and
    integral. Any feasible solution with value $s$ gives the bound
    $\Edim{\Sigma_k} \ge r - 1 + s$.
    
    Using the two above constraints we rewrite the objective function
    \[
        r + 2k -  \min \left( \sum_{j=c+1}^{n-2} (2j-n+3) \mu_j + n \mu_{n-1} \right)
    \]
    Note that all coefficients are non-negative and thus the minimum is at least
    $0$.

    For $0 < k \le r\lfloor\frac{n-2}{2}\rfloor$ set $\ell = \lceil
    \frac{k}{r} \rceil \le c+1$. Define $\mu = (\mu_0,\dots,\mu_{n-1}) \in
    \Z^{n}$ by
    \[
        \left(\begin{array}{l} \mu_{\ell-1} \\ \mu_{\ell} \end{array}\right) = 
            \begin{pmatrix} 1 & 1 \\ \ell-1 & \ell \end{pmatrix}^{-1}
            \begin{pmatrix} r \\ k \end{pmatrix} = 
            \begin{pmatrix} r\ell - k  \\ k - r( \ell-1)  \end{pmatrix} 
    \]
    and $\mu_j = 0$ otherwise. For $n$ odd we have $\ell \le
    \lfloor\frac{n-2}{2}\rfloor = c$ and $\mu$ gives a feasible solution with
    value $0$ in the minimization above.  If $n$ is even and $\ell = c+1$ the
    feasible solution yields a total value of $r + 2k - (k + r\ell) = k +
    \frac{1}{2}rn$. Note that the second case is vacuous for $n$ odd.

    For $r\lfloor\frac{n-2}{2}\rfloor < k$, let $h = k - r \lfloor
    \frac{n-2}{2} \rfloor - \alpha\lfloor\frac{n+1}{2}\rfloor$ set
    \begin{align*}
            \mu_{n-1} &= \alpha  &
            \mu_{c\phantom{+1}}   &= r - \alpha - 1 &
            \mu_{c+h\phantom{+1}}   &= 1 
    \intertext{ for $n$ odd and }
            \mu_{n-1} &= \alpha  &
            \mu_{c+1}   &= r - \alpha - 1 &
            \mu_{c+h+1}   &= 1 
    \end{align*}
    for $n$ even and $\mu_j = 0$ for all other $j$.
\end{proof}

As can be seen in the proof, the feasible solution for $\ell \le
\lfloor\frac{n-3}{2}\rfloor$ is given by a basic solution to the linear
program relaxation and it can be checked that this indeed gives the optimal
solution.  However, the coefficient for $\mu_{n-1}$ keeps this circumstance
from being true for $\ell > \lfloor\frac{n-3}{2}\rfloor$. 

In conjunction with Corollary~\ref{cor:proj_obstruct} this gives the following
definitive result concerning the non-projectability of skeleta of $\Psimplex$.

\begin{thm}
  \label{thm:ProjectionProdSimplices}
  Let $n \ge 2$ and $r \ge 1$ and set 
  $\alpha = \left\lfloor \frac{ k - r \lfloor\frac{n-2}{2}\rfloor
  }{\lfloor\frac{n+1}{2}\rfloor} \right\rfloor$. 
  If
  \[
  e <
  \left\{
  \begin{array}{l@{\text{ \ \ for \ \ }}ccccc}
      r + 2k + 1, & 
            0 & \le & k & \le & r \lfloor\tfrac{n-3}{2}\rfloor\\
      \tfrac{1}{2} r(n-2) + k + 1, &
            r\lfloor\tfrac{n-3}{2}\rfloor & < & k & \le &
            r\lfloor\tfrac{n-2}{2}\rfloor\\ 
      r(n-2)+\alpha+1, & 
            r\lfloor\tfrac{n-2}{2} \rfloor & < & k & < & r(n-1)
    \end{array}
  \right.
  \]
  then there exists no realization of the $r$-fold
  product $\Psimplex$ of $(n-1)$-simplices such that a projection to $\R^e$
  retains the $k$-skeleton.\qed
\end{thm}

For $r=1$ Theorem~\ref{thm:ProjectionProdSimplices} states that there is no
affine projection of the $(2k+2)$-simplex to $\R^{(2k+1)}$ which preserves the
$k$-skeleton. This is exactly the linear Van~Kampen--Flores Theorem. Thus, in
some sense Theorem~\ref{thm:ProjectionProdSimplices} is a generalization of
the Van~Kampen--Flores Theorem from simplices to products of simplices. As a
special case it gives yet another proof that no product of two triangles maps
linearly to a $9$-gon.

Again, let us view the statement of Theorem~\ref{thm:ProjectionProdSimplices}
in comparison with upper bounds on the embeddability dimension of the
complexes $\Skel{k}{\Psimplex}$.

\begin{prop} \label{lem:IndexProdSimplicesUpper}
    Let $\Sigma_k = \Skel{k}{\Psimplex}$ be the $k$-th coskeleton complex of the
    $r$-fold product of $(n-1)$-simplices with $n \ge 2$ and $0 \le k <
    r(n-1)$.  Then
    \[
        \Edim{ \Sigma_k } \le \min \{ 2k + 2r - 1, rn - 1 \}.
    \]
\end{prop}

\begin{proof}
   We work along the same lines as in the proof of Proposition
   \ref{prop:prod_poly_upper} and we use the fact that
   \[
        \Skel{\ell}{\simplex{n-1}} \cong \tbinom{[n]}{\le \ell+1}
        \hookrightarrow \partial\simplex{n}
    \]
    for all $0 \le \ell \le n-1$. Therefrom it follows that $\Sigma_k
    \hookrightarrow \partial\simplex{n}^{\oplus r} = \partial
    (\simplex{n^r})^{\triangle}$ and thus $\Edim{\Sigma_k} \le rn -1$.
    However, since $\dim \Sigma_k = r+k-1$ the bound given by Proposition
    \ref{prop:edim_triv_bnds} is better for $k \le \frac{1}{2} r(n-2)$.
\end{proof}

Combining the upper bounds with the lower bounds from Theorem
\ref{thm:IndexProdSimplices} yields that the result of Theorem
\ref{thm:ProjectionProdSimplices} is sharp for $k \le r
\lfloor\frac{n-3}{2}\rfloor$. On the geometric side, this is complemented in
the work of Matschke, Pfeifle, and Pilaud \cite{mpp09} on
\emph{prodsimplicial-neighborly polytopes}. The constructions given in
\cite{mpp09} yield products of simplices for which the projections retain the
$k$-skeleta for $k \le r \lfloor\frac{n-3}{2}\rfloor$.  Their constructions
also include products of simplices of different dimensions and they generalize
the topological obstructions to give bounds in the mixed case.

\section{Wedge products}
\label{sec:WedgeProducts}
\newcommand\1{\mathbf{1}}

The \emph{wedge product} $P \wedgeproduct Q$ of two polytopes $P$ and $Q$ is a
geometric degeneration of the product $Q^m$ that bears very interesting
combinatorial properties. It corresponds to an iterated subdirect product in the
sense of McMullen~\cite{mcm76} and is dual to a \emph{wreath product}
as studied by Joswig \& Lutz \cite{lutz05:_one}.

Our motivation for studying wedge products stems from the work of R\"orig \&
Ziegler \cite[Ch.~4]{Z109} on questions concerning the realizability of equivelar
surfaces. In short, an equivelar surface is a $2$-dimensional polytopal
surface that satisfies certain regularity conditions. It is both combinatorially
and geometrically challenging to construct equivelar surfaces as they exhibit
extremal combinatorial behavior. For example, unlike triangulated surfaces,
equivelar surfaces need not posses a geometric realization with flat and
convex faces (cf.~Betke~\&~Gritzmann~\cite{betke82}).

In \cite{Z109} it is shown that a certain family of wedge products $\wp{r}{n-1}$
contains equivelar surfaces in their $2$-skeleta. Furthermore, for all $r \ge 3$
the surface contained in $\wp{r}{1}$ posses a straight-line realization in
$3$-space. The approach is to give a geometric realization of $\wp{r}{1}$ such
that a projection to $\R^4$ strictly preserves the surface and the resulting
polytope carries the surface in its lower hull.

In this section we prove non-projectability results regarding skeleta of wedge
products and, in particular, we show that for $r \ge 4$ and $n \ge 3$ there is
no realization of the wedge product $\wp{r}{n-1}$ such that a projection to
$\R^4$ strictly preserves the equivelar surface. 

\subsection{Wedge products and products}
\label{ssec:wedge-prod-prod}

Wedge products of polytopes were introduced in \cite{Z109} from several
perspectives such as an iteration of a generalized wedge construction and in
terms of interior and exterior presentations. In this paper, we will only need
the description in terms of facet-defining halfspaces.

\begin{dfn}[{\cite[Def.~4.10]{Z109}}]
For polytopes $P = \{ y \in \R^d : a_i \cdot y \le 1 \text{ for all } i =
1,\dots,m \}$ and $Q = \{ x \in \R^{d^\prime} : B x \le \1 \}$
the \emph{wedge product} of $P$ and $Q$ is the polytope
\[
    P \wedgeproduct Q :=
    \left\{ (x_1,\dots,x_m,y) \in (\R^{d^\prime})^m \times \R^d : B x_i \le (1 -
    a_i \cdot y) \1 \text{ for all } i = 1,\dots, m \right\}.
\]
\end{dfn}

The geometry and the combinatorics of wedge products are studied in \cite{Z109}.
For our purposes it is sufficient to know the combinatorial type of $P
\wedgeproduct Q$ in the form of intersections of facets.

\begin{thm}\label{thm:wp_lattice}
    Let $P$ and $Q$ be polytopes with facets indexed by $[m]$ and $[n]$,
    respectively. The face lattice of $P \wedgeproduct Q$ is given by the
    collection of tuples $(H_1,\dots,H_m)$ with $H_1,\dots,H_m \subseteq [n]$
    such that
    \begin{enumerate}
        \item $H_i = I_Q(F_i)$ for some face $F_i \subseteq Q$ for all $i$,
            and
        \item $\{ j \in [m] : H_j = [n] \} = I_P(G)$ for some face $G
            \subseteq P$.
    \end{enumerate}
    The order relation is given by componentwise reverse inclusion. The
    dimension of the face $(H_1,\dots,H_m)$ is given by $\sum \dim F_i + \dim
    G$. 
\end{thm}
\begin{proof}
    It follows from the lattice structure of $\fl(P)$ and $\fl(Q)$ that the stated
    poset is a atomic and coatomic lattice.  It is known that two
    atomic-coatomic lattices are isomorphic if and only if they have isomorphic
    atom-coatom incidences.  The bijection on the collection of facets is clear
    and the vertices are determined by Theorem~4.13 in \cite{Z109} and
    correspond to admissible tuples $(H_1,H_2,\dots,H_m)$ with $F_i$ of
    dimension at most $0$ and $G$ a vertex.
\end{proof}

An alternative approach to wedge products and Theorem~\ref{thm:wp_lattice}
appears in \cite[Thm.~2.21]{raman-thesis}. The following observation links the wedge product to the usual product.

\begin{prop}[{\cite[Prop.~4.12]{Z109}}]
    The intersection of the wedge product $P \wedgeproduct Q$ with the linear
    space $L = (\R^{d^\prime})^m \times \{ 0 \} \subset  (\R^{d^\prime})^m
    \times \R^d$ is affinely isomorphic to $Q^m$. In particular, the
    intersection is given by the faces $(H_1,\dots,H_m)$ with $H_i \not= [n]$.
\end{prop}

It follows from Proposition~\ref{thm:wp_lattice} that every $k$-face for
$k \ge 0$ of the product $Q^m$ is the unique intersection of $L$ with a face
of dimension $k+\dim P$ of $P \wedgeproduct Q$. 

\begin{lem}
    \label{lem:WedgeProdSkeleta}
    Let $P$ and $Q$ be polytopes with $m$ being the number of facets of $P$.
    Then for any $0 \le k \le m\,\dim Q$ we have
    \[
    \Skel{k}{Q^m} \hookrightarrow \Skel{k + \dim P}{P \wedgeproduct Q}.
    \]
 \qed
\end{lem}

We call the image of the $k$-skeleton of the product $Q^m$ in $P \wedgeproduct
Q$ the \emph{special} $(k+\dim P)$-faces of the wedge product. These special
faces cover all vertices of the wedge product by Theorem~\ref{thm:wp_lattice}. 
The bottom line is that we can re-use the bounds obtained in
Section~\ref{sec:ProductsOfSimplices} to handle projections of wedge products of
polygons and simplices.

\subsection{Projections of wedge products of polygons and simplices}
\label{ssec:projectability}

In the following we restrict ourselves to the wedge product $\wp{r}{n-1} = \P_r
\wedgeproduct \simplex{n-1}$ of an $r$-gon $\P_r$ and an $(n-1)$-simplex
$\simplex{n-1}$. It follows from Theorem~\ref{thm:wp_lattice} that $\wp{r}{n-1}$
is an $(r(n-1)+2)$-dimensional polytope with $rn$ facets. Using the correspondence
established in Lemma~\ref{lem:WedgeProdSkeleta} we apply the result of
Section~\ref{sec:ProductsOfSimplices} to the projectability of the $k$-skeleta
of wedge products.

\begin{prop}
  \label{prop:ObstructionSpecialKFaces}
  There exists no realization of the wedge product $\wp{r}{n-1}$ of $r$-gon and
  $(n-1)$-simplex with $r\ge 4$ and $n \ge 2$ such that the projection to $\R^e$
  preserves its special $k$-faces for $k\ge 2$ if
  \[
  e <
  \left\{
    \begin{array}{llccccc}
      r+2k-1 &\text{if}&2&\le& k &\le&r\lfloor\tfrac{n-3}{2}\rfloor + 2\\
      \tfrac{1}{2} r(n-2) + k + 1&\text{if}
      &r\lfloor\tfrac{n-3}{2}\rfloor + 2 & < & k
      &\le&r\lfloor\tfrac{n-2}{2}\rfloor + 2\\ 
      r(n-2) + \alpha + 3&\text{if}&r\lfloor\tfrac{n-2}{2}\rfloor + 2 &<&
      k&<&r(n-1)+2 \\
    \end{array}
  \right.
  \]
  and
  \[
  \alpha = 
  \left\lfloor
    \frac{k-2-\left\lfloor\frac{n-2}{2}\right\rfloor}
    {\left\lfloor\frac{n+1}{2}\right\rfloor}
  \right\rfloor
  .
  \]
\end{prop}
\begin{proof}
    We are able to apply Theorem~\ref{thm:main} since the special faces cover
    all vertices and every face of a strictly preserved face is also strictly
    preserved.
    The strict projection complex of a projection strictly preserving the special
    $k$-faces of the wedge product contains the $(k-2)$-nd coskeleton complex of
    the product $\Psimplex$ (see Lemma~\ref{lem:WedgeProdSkeleta}). Hence the
    embeddability dimension of the special $k$-faces of the wedge product is
    equal to the embeddability dimension of $\Skel{k-2}{\Psimplex}$ given by
    Theorem~\ref{thm:IndexProdSimplices}. Plugging these bounds into
    Corollary~\ref{cor:proj_obstruct} we obtain:
  \[
  e 
  <   \Edim{\Skel{k-2}{\Psimplex}} - r + 4 
  \le \Edim{\Skel{k}{\wp{r}{n-1}}} + (r(n-1) + 2) - rn + 2.\qedhere
  \]
\end{proof}

Since the $k$-skeleton of the wedge product obviously contains the special
$k$-faces we obtain the following result for the projectability of skeleta of
the wedge product $\wp{r}{n-1}$.

\begin{thm}
  \label{thm:ObstructionWedgeProductSkeleta}
  There exists no realization of the wedge product $\wp{r}{n-1}$ of $r$-gon and
  $(n-1)$-simplex with $r\ge 4$ and $n \ge 2$ such that the projection to $\R^e$
  preserves the $k$-skeleton for $k \ge 0$ if
  \[
  e <
  \left\{
    \begin{array}{llccccc}
      r+2k-1 &\text{if}&2&\le& k &\le&r\lfloor\tfrac{n-3}{2}\rfloor + 2\\
      \tfrac{1}{2} r(n-2) + k + 1&\text{if}
      &r\lfloor\tfrac{n-3}{2}\rfloor + 2 & < & k
      &\le&r\lfloor\tfrac{n-2}{2}\rfloor + 2\\ 
      r(n-2) + \alpha + 3&\text{if}&r\lfloor\tfrac{n-2}{2}\rfloor + 2 &<&
      k&<&r(n-1)+2 \\
    \end{array}
  \right.
  \]
  and
  \[
  \alpha = 
  \left\lfloor
    \frac{k-2-\left\lfloor\frac{n-2}{2}\right\rfloor}
    {\left\lfloor\frac{n+1}{2}\right\rfloor}
  \right\rfloor
  .
  \]
\end{thm}
\begin{proof}
    The vertices of the wedge product correspond to the vectors $\H_V$
    given by:
    \begin{equation}
        \label{eq:WPVerts}
        \H_V = \left\{ (H_1,\dots,H_r) \in \wp{r}{n-1} : H_i
          \not= [n] \Rightarrow |H_i| = n-1 \right\}.
    \end{equation}
    We pick a subfamily of vertices corresponding to the vectors
    $([n],[n],H_3,\dots,H_r)$ with $|H_i|=n-1$ for $i = 3,\dots,r$. Considering
    only the last $r-2$ components of the vector we obtain the following
    inclusion of coskeleton complexes:
    \[
    \Skel{0}{\simplex{n-1}^{r-2}} \hookrightarrow \Skel{0}{\wp{r}{n-1}}.
    \]
    The embeddability dimension of $\Skel{0}{\simplex{n-1}^{r-2}}$ is $2r-5$ by
    Theorem~\ref{thm:IndexProdSimplices}. Thus we obtain the following bound
    on $e$ with Corollary~\ref{cor:proj_obstruct}:
    \[
    \Edim{\Skel{0}{\wp{r}{n-1}}} + r(n-1)+2 - rn + 2 \ge
    \Edim{\Skel{0}{\simplex{n-1}^{r-2}}} - r + 4 = 
    r - 1 > e.
    \]
    The $1$-skeleton of the wedge product contains a subfamily of edges
    corresponding to the vectors $([n],H_2,\dots,H_r)$. As for the vertices we
    obtain an inclusion of coskeleton complexes:
    \[
    \Skel{0}{\simplex{n-1}^{r-1}} \hookrightarrow \Skel{1}{\wp{r}{n-1}}.
    \]
    Since by Theorem~\ref{thm:IndexProdSimplices} the embeddability dimension of
    $\Skel{0}{\simplex{n-1}^{r-1}}$ is $2r-3$ we obtain the following bound for
    the dimension projected onto using Corollary~\ref{cor:proj_obstruct}:
    \[
    \Edim{\Skel{1}{\wp{r}{n-1}}} - r + 4 \ge
    \Edim{\Skel{0}{\simplex{n-1}^{r-1}}} - r + 4 = 
    r + 1 > e.
    \]
    For $k \ge 2$ we simply use Proposition~\ref{prop:ObstructionSpecialKFaces}.
\end{proof}

\subsection{Equivelar surfaces in wedge products}

It is shown in \cite{Z109} that the wedge product $\wp{r}{n-1} = \P_r
\wedgeproduct \simplex{n-1}$ of an $r$-gon and an $(n-1)$-simplex carries a
very interesting equivelar surface $\WPsurf{r}{2n}$ in its $2$-skeleton. A
main result of \cite{Z109} is that in some cases this combinatorial embedding
can be used to obtain a geometric embedding in $3$-space.
Using the machinery developed in Section~\ref{sec:projections} and the results
of Section~\ref{sec:ProductsOfSimplices} we complement the above result about
projections of equivelar surfaces.

The $2$-skeleton of the wedge product $\wp{r}{n-1}$ is a
furtile ground for embedding equivelar surfaces. Consider the special $2$-faces of Lemma~\ref{lem:WedgeProdSkeleta}. They correspond to:
\[
    \H_R =
    \left\{
        (H_1,\dots,H_r) \in \wp{r}{n-1} :
        |H_i| = n-1 \text{ for all } i \in [r] 
    \right\}.
\]
So for every choice $j_1,\dots,j_r \in [n]$ the tuple
\[
    H = ([n] \setminus j_1, [n] \setminus j_2, \dots, [n] \setminus
    j_r)
\]
represents a special $2$-face of $\wp{r}{n-1}$ and each such face is isomorphic
to an $r$-gon. Indeed, for every $i \in [r]$ the tuple
\[
    H^i = ([n] \setminus j_1, \dots, [n] \setminus j_{i-1}, [n], 
    [n] \setminus j_{i+1}, \dots, [n] \setminus j_r)
\]
corresponds to an edge of $H$ by Theorem~\ref{thm:wp_lattice} and hence $H$ is a
$2$-dimensional face with $r$ edges. We denote the collection of these $r$-gon
edges by $\H_E$: They correspond to tuples $(H_1,\dots,H_r)$ with $|H_i| = n -1$
for all but a unique $i_0 \in [r]$ with $H_{i_0} = [n]$.

In~\cite{Z109} the following subcomplex of the wedge product $\wp{r}{n-1}$ is
discussed: For $r \ge 3$ and $n \ge 2$ consider the subcomplex $\WPsurf{r}{2n}$
generated by the following collection of $r$-gons of the wedge product
$\wp{r}{n-1}$:
\[
    \WPsurf{r}{2n} =
    \Big\{ ([n] \setminus j_1,\ldots,[n]\setminus j_{r}) :
    \sum_{k=1}^{r} j_k \equiv 0,1 \mod n\Big\} \subseteq
    \H_R.
\]
The subcomplex $\WPsurf{r}{2n}$ contains all the vertices
and all the edges of $\H_E$ of $\wp{r}{n-1}$. It is shown in~\cite{Z109} that $\WPsurf{r}{2n}$ is
a regular (polyhedral) surface~$\WPsurf{r}{2n}$ of type $\{r,2n\}$, i.e.\ an
(orientable) polyhedral $2$-manifold that is
\begin{compactitem}
\item \emph{equivelar}: all faces are $r$-gons and every vertex is incident to
    $2n$ faces, and even 
\item \emph{regular}: the automorphism group acts transitively on the flags of
    the surface.
\end{compactitem}

For the special case $n=2$ there are deformed realizations of the wedge products
$\wp{r}{1}$ and projections that yield embeddings of the surfaces
$\WPsurf{r}{4}$ in $\R^3$.

\begin{thm}[{\cite[Thm.~4.26]{Z109}}]
  The wedge product $\wp{r}{1}$ has a realization such that all the faces
  corresponding to the surface $\WPsurf{r}{4}\subset \wp{r}{1}$ are preserved by
  the projection to $\R^4$. Hence there is a realization of $\WPsurf{r}{4}$ in $\R^3$.
\end{thm}

So there was hope that some realizations of the wedge products for other
parameters $r$ and $n$ would yield realizations in $\R^3$ as well. But with the
techniques developed in this article we obtain the following negative result.

\begin{thm}
  \label{thm:ObstructionWedgeProductSurface}
  There is no realization of the wedge product $\wp{r}{n-1}$, with $n \ge 3$ and
  $r \geq 4$, such that all the faces corresponding to the surface
  $\WPsurf{r}{2n}$ are strictly preserved by the projection
  $\pi:\R^{2+r(n-1)}\rightarrow\R^{e}$ for $e < r+1$.
\end{thm}
\begin{proof}
    We prove the theorem by contradiction. So assume that there exists a
    realization of $\wp{r}{n-1}$ such that the surface $\WPsurf{r}{2n}$ is
    strictly preserved by the projection to $\R^e$ with $e < r+1$. By
    Theorem~\ref{thm:main} the embeddability dimension of the strict projection
    complex $\K = \K(\wp{r}{n-1},\pi)$ is then
    \begin{equation}
        \Edim{\K} \le rn-(r(n-1)+2)+e-2 = r + e - 4 < 2r-3.
        \tag{\sf WP}
        \label{eq:edim_wp}
    \end{equation}
    Since the polygons of the wedge product surface $\WPsurf{r}{2n}$ are
    strictly preserved by the projection $\pi$ the simplicial complex $\K$
    contains a subcomplex $\Sigma$ corresponding to the polygons of
    $\WPsurf{r}{2n}$. The strict projection complex of all special $r$-gons is
    $\Skel{0}{\Psimplex}$ by Lemma~\ref{lem:WedgeProdSkeleta}. Hence
    \[
    \Sigma = \{(j_1,\dots,j_r)\,|\, \sum_{k=1}^{r} j_k \equiv 0,1 \mod n\}
    \ \subset\ \Skel{0}{\Psimplex}.
    \]
    We remove the asymmetry from $\Sigma$ by only considering the edges
    $([n],[n]\setminus j_2,\dots,[n]\setminus j_{r})$ for $j_i \in [n]$ of the
    wedge product. The strict projection complex of these edges is
    $\Skel{0}{\simplex{n-1}^{r-1}}$. By Theorem~\ref{thm:IndexProdSimplices} the
    embeddability dimension of $\Skel{0}{\simplex{n-1}^{r-1}}$ is $2r-3$. Hence
    the embeddability dimension of $\K$ is at least $2r-3$
    because $\Skel{0}{\simplex{n-1}^{r-1}} \subseteq \Sigma \subseteq
    \K(\wp{r}{n-1},\pi)$. This is a contradiction to
    Equation~\eqref{eq:edim_wp}. So there exists no realization of $\wp{r}{n-1}$
    such that the surface $\WPsurf{r}{2n}$ is strictly preserved by the
    projection to $\R^e$.
\end{proof}

Theorem~\ref{thm:ObstructionWedgeProductSurface} does not obstruct
straight-line realizations of the surfaces $\WPsurf{r}{2n}$ in $\R^3$ in
general. It only exhibits the limitations of the approach taken in \cite{Z109}.

\bibliographystyle{siam}
\begin{small}
  \bibliography{PolytopeSkeleta}
\end{small}

\end{document}